\documentclass{amsproc}
\usepackage{epsfig}
\newtheorem{theorem}{Theorem}[section]

\theoremstyle{definition}

\newtheorem{corollary}[theorem]{Corollary}

\theoremstyle{remark}
\newtheorem{remark}[theorem]{Remark}

\numberwithin{equation}{section}

\begin{document}

\title{An explicit description of the simplicial group $K(A,n)$}
\author{Mihai D. Staic}
\address{Department of Mathematics, DePaul University,  Chicago, IL 60614, USA }
\address{Institute of Mathematics of the Romanian Academy, PO.BOX 1-764, RO-70700 Bu\-cha\-rest, Romania.}
\thanks{Research partially supported by the
CNCSIS project ''Hopf algebras, cyclic homology and monoidal categories'' contract nr. 560/2009, CNCSIS code $ID_{-}69$.}
\email{mstaic@depaul.edu}



\date{January 1, 1994 and, in revised form, June 22, 1994.}



\begin{abstract} We give a new explicit construction for the simplicial group  $K(A,n)$ and explain its topological interpretation.

\end{abstract}

\maketitle


%

\section*{Introduction}

Simplicial groups are purely algebraic objects that are used in algebraic to\-po\-lo\-gy to formulate classification 
results. Just like for topological spaces one can talk about the n-th homotopy group of a simplicial group. A $K(A,n)$ simplicial group is determined by the fact that  $\pi_i(K(A,n))=0$ if $i\neq n$ and $\pi_n(K(A,n))=A$.  In other words, it is the algebraic object  corresponding  to an  Eilenberg-MacLane space $K(A,n)$. If $A$ is a fixed commutative group there is an iterating procedure that gives a  presentation of the simplicial group of $K(A,n)$ (see  \cite{may}). Unfortunately  some of the topological nature of simplicial objects is lost in the iterating
process.  There are also explicit description of $K(A,n)$ (see for example \cite{gj} or \cite{w}), but again the topological flavor is not transparent.

In this paper we give a new explicit description of the simplicial group $K(A,n)$. The main advantage of our presentation is that it has a nice topological interpretation. Also the description is very simple and is presented in terms of the generating map of the simplicial category $\Delta$.

In the first section we recall basic  definitions, properties and examples of simplicial groups. The second section starts with the description of $K(A,2)$. This construction appears  in non explicit way in \cite{sm1} and it was the starting point of the paper. We show that $K(A,2)$ is  a cyclic object. For a better understanding of our general  construction we also present the case of $K(A,3)$.  In the third section  we give the description of the simplicial group $K(A,n)$. The punch line is in the way we chose the index for the elements of the group $K(A,n)_q=A^{\,(^q_n)}$. We explain the topological interpretation  and discuss some possible applications.

The last section deals with a similar construction in the context of  Hopf Algebras. More exactly, to every commutative Hopf algebra $H$ we associate a cyclic object $\,_2K(H)$. If $H$ is the group algebra $k[A]$ associated to a commutative group $A$, then the cyclic object  $\,_2K(H)$ is just the linearization of the cyclic object $K(A,2)$ mentioned above.

\section{Preliminaries}

We recall from \cite{may} and \cite{lo} a few facts about simplicial groups.
First the definition of the simplicial category $\Delta$. The objects in $\Delta$ are the finite ordered sets $\overline{n}=\{0,1,...,n\}$. The morphisms are the order preserving maps. One can show the any morphism in $\Delta$ can be written as a composition of maps $d^i:\overline{n} \to \overline{n+1}$ and $s^i:\overline{n} \to \overline{n-1}$ where
$$d^i(u)= \left\{
\begin{array}{c l}
u\; &{\rm if}\;  u< i\\
u+1 &{\rm if}\;  u\geq i
\end{array}
\right.
$$
$$s^i(u)= \left\{
\begin{array}{c l}
u\; &{\rm if}\;  u\leq i\\
u-1 &{\rm if}\;  u> i
\end{array}
\right.
$$
A simplicial group is a functor $K:\Delta^{op}\to Gr$. More explicitly, a simplicial group is a set of groups $X_q$, $q\geq 0$ together with a colection of group  morphisms $\partial_i:X_q\to X_{q-1}$ and $s_i:X_q\to X_{q+1}$ for all $0\leq i\leq q$ such that the following identities hold
\begin{eqnarray*}
&&\partial_i\partial_j=\partial_{j-1}\partial_i\; {\rm if}\; i<j\\
&&s_is_j=s_{j+1}s_i\; {\rm if}\; i\leq j\\
&&\partial_i s_j=s_{j-1}\partial_i\; {\rm if}\; i<j\\
&&\partial_js_j=\partial_{j+1}s_j=id\\
&&\partial_i s_j=s_j\partial_{i-1}\; {\rm if}\; i>j+1\\
\end{eqnarray*}
Let ${\bf K}=(K_q,\partial_i,s_i)$ be a simplicial group.  We denote
$$\overline{K}_q=K_q\cap Ker(\partial_0) \cap ... \cap Ker(\partial_{q-1}).$$
One can show that $\partial_{q+1}(\overline{K}_{q+1})\subseteq \overline{K}_{q+1}$ and so $\overline{K}$ is a chain complex. The homotopy groups of the simplicial group ${\bf K}$ are defined by
$$\pi_q({\bf K})=H_q(\overline{\bf K}).$$
If $q\geq 2$ then $\pi_q({\bf K})$ is an abelian group. If ${\bf K}$ has the property that $\pi_i({\bf K})=0$ if $i\neq n$ and
$\pi_n({\bf K})=G$ then it is called a Eilenberg-MacLane simplicial group and is denoted by  $K(G,n)$.

We recall  from \cite{may} the construction of the $K(G,1)$ simplicial group. Define
$$K_q=G^q$$
where the elements of $G^q$ are q-tuples $(g_0,g_1,...,g_{q-1})$.
For every $0\leq i\leq q-1$ define
$$\partial_i:G^q\to G^{q-1}\; \; {\rm and}\; \;  s_i:G^q\to G^{q+1}$$
\begin{eqnarray*}
\partial_0(g_0,g_1,...,g_{q-1})&=&(g_1,...,g_{q-1})\\
\partial_i(g_0,g_1,...,g_{q-1})&=&(g_0,..,g_{i-1}g_i,...,g_{q-1})\\
\partial_q(g_0,g_1,...,g_{q-1})&=&(g_0,...,g_{q-2})
\end{eqnarray*}
and
\begin{eqnarray*}
s_0(g_0,g_1,...,g_{q-1})&=&(1,g_0,...,g_{q-1})\\
s_i(g_0,g_1,...,g_{q-1})&=&(g_0,..,g_{i-1},1,g_{i},...,g_{q-1})\\
s_q(g_0,g_1,...,g_{q-1})&=&(g_0,...,g_{q-1},1)
\end{eqnarray*}
One can show that $(K_q,\partial_i,s_i)$ is a $K(G,1)$ simplicial group. Moreover if 
$$\tau_q:K_q=G^q\to K_{q}=G^{q}$$
\begin{eqnarray*}
\tau_q(g_0,g_1,...,g_{q-1})=((g_0g_1...g_{q-1})^{-1},g_0,...,g_{q-2})
\end{eqnarray*}
then $(K_q,\partial_i,s_i, \tau_q)$ is a cyclic simplicial group.

Let $(H,\Delta, \varepsilon)$ be a Hopf algebra. We use the Swedler's sigma notation $\Delta(h)=h^{<1>}\otimes h^{<2>}$. A pair $(\delta,\sigma)$ is a modular pair in involution if  $\sigma\in H$ is group-like element,  $\delta:H\to k$  is a character for $H$  and  $\tilde{S}_{\sigma}^2=id$, where $\tilde{S}_{\sigma}(h)=\sigma\Sigma\delta(h^{<2>})S(h^{<1>})$.  It was proved in \cite{kr} that a Hopf algebra $H$ and any modular pair in involution $(\delta,\sigma)$ one can associate a cyclic module $H_n^{(\delta,\sigma)}$. First the simplicial structure $H_n^{(\delta,\sigma)}=H^{\otimes n}$
\begin{eqnarray*}
\delta_0(h_1\otimes h_2\otimes...\otimes h_n)&=&\varepsilon(h_1)h_2\otimes...\otimes h_n\\
\delta_i(h_1\otimes h_2\otimes...\otimes h_n)&=&h_1\otimes...\otimes h_ih_{i+1}\otimes ... \otimes h_n\\
\delta_n(h_1\otimes h_2\otimes...\otimes h_n)&=&\delta(h_n)h_1\otimes...\otimes h_{n-1}\\
\sigma_0(h_1\otimes h_2\otimes...\otimes h_n)&=&1\otimes h_1\otimes...\otimes h_n\\
\sigma_i(h_1\otimes h_2\otimes...\otimes h_n)&=&h_1\otimes...h_i\otimes 1\otimes h_{i+1}...\otimes h_n\\
\sigma_n(h_1\otimes h_2\otimes...\otimes h_n)&=&1\otimes h_1\otimes...\otimes h_n\\
\end{eqnarray*}
For the cyclic action define
\begin{eqnarray*}
 \tau_n(h_1\otimes h_2\otimes ... \otimes h_n)=\delta(h_n^{<2>})\sigma S(h_1^{<1>}h_2^{<1>}...h_n^{<1>})\otimes h_1^{<2>}\otimes ...\otimes h_{n-1}^{<2>} \label{cycl}
\end{eqnarray*}
\begin{remark}
When we specialize to the case $H=k[G]$ the group algebra associated to a group $G$, $\sigma=1$ and $\delta=\varepsilon$ we get the linearization of the cyclic simplicial group $K(G,1)$ described above.
\end{remark}
\begin{remark} If $H$ is a cocomutative Hopf algebra then $(\varepsilon,1)$ is a modular pair in involution.  The simplicial object $H^{(\varepsilon, 1)}$ has a natural structure of symmetric object (see \cite{lo} for the definition of symmetric objects). The action of the symmetric group $\Sigma_{n+1}$ is given by
\begin{eqnarray*}
(1,2)(h_1\otimes ...\otimes h_n)&=&S(h_1^{<1>})\otimes h_1^{<2>}h_2\otimes...\otimes h_n\\
(i,i+1)(h_1\otimes ...\otimes h_n)&=&h_1\otimes...\otimes h_{i-1}h_i^{<1>}\otimes S(h_i^{<2>})\otimes h_i^{<3>}h_{i+1}\otimes ... \otimes h_n\\
(n,n+1)(h_1\otimes ...\otimes h_n)&=&h_1\otimes...\otimes h_{n-1}h_n^{<1>}\otimes S(h_n^{<2>}).
\end{eqnarray*}
One can notice that the map $\tau_n$  is given by the action of the cycle $(1,2,...,n+1)\in \Sigma_{n+1}$ on $H^{\otimes n}$. And so the cyclic structure on $H^{(\varepsilon, 1)}$ is induced by the above symmetric structure.
\end{remark}

\section{Eilenberg-MacLane simplicial groups (case $n=2$ and $n=3$)}
\subsection{$K(A,2)$}
In this section $A$ is a commutative group. In \cite{sm1} we introduced the secondary cohomology $\,_2H^n(A,B)$ (where $B$ is  a commutative group).  We proved that to a topological space $X$ with $\pi_1(X)=1$ one can associate an invariant $\,_2\kappa^4\in  \,_2H^4(\pi_2(X),\pi_3(X))$. As a byproduct of that that construction one gets an explicit description of simplicial group $K(A,2)$. The basic idea is that in order to color the 2-dimensional skeleton of the q-simplex $\Delta_q$ with elements  of the group $A$, it is enough to say what the colors are for all 2-simplices of the form $[u,v,v+1]$ where $0\leq u<v\leq q-1$. For the rest of the 2-faces the color is determined by ''homotopy''.  Take for example $\Delta_3$, once we fix the color of $[0,1,2]$, $[0,2,3]$ and $[1,2,3]$ as $a_{0,1}$, $a_{0,2}$ and $a_{1,2}$ respectively, then the face $[0,1,3]$ must have the color $a_{0,1}a_{0,2}a_{1,2}^{-1}$. To see this we should think that $A$ is the second homotopy group $\pi_2(X)$ of a simply connected topological space $X$. We send the 1-skeleton of $\Delta_3$ in the fixed point of $X$, and three of the 2-faces of $\Delta_3$ according to the above prescription. If we want to have a map from $\Delta_3$ to $X$ we are forced to send the face $[0,1,3]$ to the element $a_{0,1}a_{0,2}a_{1,2}^{-1}$.

Define
$$K_q=A^{\frac{q(q-1)}{2}}$$
The elements of  $A^{\frac{q(q-1)}{2}}$ are $\frac{q(q-1)}{2}$-tuples $(a_{u,v})_{(0\leq u<v\leq q-1)}$ with the index in the lexicographic order:
$$(a_{0,1},a_{0,2},...,a_{0,q-1},a_{1,2},a_{1,3},...,a_{1,q-1},...,a_{q-2,q-1})$$
For every $0\leq i\leq q$ we define
$$\partial_i:K_q=A^{\frac{q(q-1)}{2}}\to K_{q-1}=A^{\frac{(q-1)(q-2)}{2}}$$
$\partial_i((a_{u,v})_{(0\leq u<v\leq q-1)})=(b_{u,v})_{(0\leq u<v\leq q-2)}$ where
$$b_{u,v}= \left\{
\begin{array}{c l}
a_{u,v}\; &{\rm if}\; 0\leq u<v < i-1\\
a_{u,v}a_{u,i}a_{v,i}^{-1}\; &{\rm if}\; 0\leq u<v= i-1\\
a_{u,v+1}\; &{\rm if}\; 0\leq u\leq i-1 < v\\
a_{u+1,v+1}\; &{\rm if}\; i-1<u <v
\end{array}
\right.
$$
and
$$s_i: K_q=A^{\frac{q(q-1)}{2}}\to K_{q+1}=A^{\frac{(q+1)q}{2}}$$
$s_i((a_{u,v})_{(0\leq u<v\leq q-1)})=(c_{u,v})_{(0\leq u<v\leq q)}$ where

$$c_{u,v}= \left\{
\begin{array}{c l}
a_{u,v}\; &{\rm if}\; 0\leq u<v < i\\
1\; &{\rm if}\; 0\leq u<v = i\\
a_{u,v-1}\; &{\rm if}\; 0\leq u<i<v\\
a_{u,v-1}\; &{\rm if}\; 0\leq u=i< v-1\\
1\; &{\rm if}\; 0\leq u=i=v-1\\
a_{u-1,v-1}\; &{\rm if}\; 0\leq i<u <v
\end{array}
\right.
$$

A long but straightforward computation shows that $(K_q,\partial_i, s_i)$ is a simplicial
group. Moreover one can see that
$$\overline{K}_q=K_q\cap Ker(\partial_0) \cap ... \cap Ker(\partial_{q-1})= \left\{
\begin{array}{c l}
1\; {\rm if}\; q \neq 2\\
A\; {\rm if}\; q=2
\end{array}
\right.
$$
We have
\begin{theorem} $(K_q,\partial_i, s_i)$ is a $K(A,2)$ simplicial group.
\end{theorem}
\begin{corollary} $(K_q,\partial_i, s_i,\tau_q)$ is a cyclic simplicial group, where $\tau_q:K_q\to K_q$
$$\tau_q((a_{u,v})_{(0\leq u<v\leq q-1)})=(e_{u,v})_{(0\leq u<v\leq q-1)}$$
$$e_{u,v}= \left\{
\begin{array}{c l}
a_{v-1,v}a_{v-1,v+1}...a_{v-1,q-1}a_{v,v+1}^{-1}a_{v,v+2}^{-1}...a_{v,q-1}^{-1}\; {\rm if}\; 0=u<v\\
a_{u-1,v-1}\; {\rm if}\; 0< u<v\\
\end{array}
\right.$$
\end{corollary}

\subsection{$K(A,3)$}
Just like above one can see that in order to color the 3-skeleton of the q-simplex $\Delta_q$, it is enough to color all 3-simplices  of the form $[u,v,w,w+1]$, where $0\leq u<v<w\leq q-1$. The color for the other 3-faces is determined by ''homotopy''.

We define
$$K_q=A^{\frac{q(q-1)(q-2)}{6}}$$
The elements of  $A^{\frac{q(q-1)(q-2)}{6}}$ are $\frac{q(q-1)(q-2)}{6}$-tuples $(a_{u,v,w})_{(0\leq u<v<w\leq q-1)}$ with the index in the lexicographic order:
$$(a_{0,1,2},a_{0,1,3},...,a_{0,1,q-1},a_{0,2,3},a_{0,2,4},...,a_{q-3,q-2,q-1})$$
For every $0\leq i\leq q$ we put
$$\partial_i:K_q=A^{\frac{q(q-1)(q-2)}{6}}\to K_{q-1}=A^{\frac{(q-1)(q-2)(q-3)}{6}}$$
$\partial_i((a_{u,v,w})_{(0\leq u<v<w\leq q-1)}=(b_{u,v,w})_{(0\leq u<v<w\leq q-2)}$ where
$$b_{u,v,w}= \left\{
\begin{array}{c l}
a_{u,v,w}\; &{\rm if}\; 0\leq u<v<w < i-1\\
a_{u,v,w}a_{u,v,i}a_{u,w,i}^{-1}a_{v,w,i}\; &{\rm if}\; 0\leq u<v<w= i-1\\
a_{u,v,w+1}\; &{\rm if}\; 0\leq u<v\leq i-1 < w\\
a_{u,v+1,w+1}\; &{\rm if}\; 0\leq u\leq i-1 <v< w\\
a_{u+1,v+1,w+1}\; &{\rm if}\; i-1<u <v<w
\end{array}
\right.
$$
and
$$s_i: K_q=A^{\frac{q(q-1)(q-2)}{6}}\to K_{q+1}=A^{\frac{(q+1)q(q-1)}{6}}$$
$s_i((a_{u,v,w})_{(0\leq u<v<w\leq q-1)})=(c_{u,v,w})_{(0\leq u<v<w\leq q)}$ where
$$c_{u,v,w}= \left\{
\begin{array}{c l}
a_{u,v,w}\; &{\rm if}\; 0\leq u<v<w < i\\
1\; &{\rm if}\; 0\leq u<v<w = i\\
a_{u,v,w-1}\; &{\rm if}\; 0\leq u<v<i<w\\
a_{u,v,w-1}\; &{\rm if}\; 0\leq u<v=i< w-1\\
1\; &{\rm if}\; 0\leq u<v=i=w-1\\
a_{u,v-1,w-1}\; &{\rm if}\; 0\leq u<i<v <w\\
a_{u,v-1,w-1}\; &{\rm if}\; 0\leq u=i<v-1 <w-1\\
1\; &{\rm if}\; 0\leq u=i=v-1<w-1\\
a_{u-1,v-1,w-1}\; &{\rm if}\; 0\leq i<u<v <w\\
\end{array}
\right.
$$
One can check that in this way we get a $K(A,3)$ simplicial group.

\section{Eilenberg-MacLane simplicial groups $K(A,n)$}
The above two examples suggest a  description for all  Eilenberg-MacLane simplicial groups $K(A,n)$. This time, to color the n-skeleton of the q-simplex $\Delta_q$, it is enough to color all n-simplices  of the form $[u_1,u_2,...,u_n,u_n+1]$ where $0\leq u_1<u_2<...<u_n\leq q-1$.

We define
$$K(n)_q=A^{(\,^q_n)}$$
The elements of  $A^{(\,^q_n)}$ are $(\,^q_n)$-tuples $(a_{u_1,...,u_n})_{(0\leq u_1<...<u_n\leq q-1)}$ with the index in the lexicographic order
$$(a_{0,1,...,n-2,n-1},a_{0,1,...,n-2,n},...,a_{q-n,...q-2,q-1})$$
For every $0\leq i\leq q$ define
$$\partial_i:K(n)_q=A^{(\,^q_n)}\to K(n)_{q-1}=A^{(\,^{q-1}_{\;\;n})}$$
$\partial_i((a_{u_1,...,u_n})_{(0\leq u_1<...<u_n\leq q-1)})=(b_{u_1,...,u_n})_{(0\leq u_1<...<u_n\leq q-2)}$ where
$$b_{u_1,...,u_n}= \left\{
\begin{array}{c l}
a_{u_1,...,u_n}a_{u_1,...,u_{n-1},i}^{(-1)^0}a_{u_1,...,u_{n-2},u_n,i}^{(-1)^1}...a_{u_2,...,u_n,i}^{(-1)^{n-1}}\; &{\rm if}\;  u_n=i-1\\
a_{d^i(u_1),...,d^i(u_n)}\; &{\rm if}\; u_n\neq i-1
\end{array}
\right.
$$
The degeneracy maps are
$$s_i: K(n)_q=A^{(\,^q_n)}\to K(n)_{q+1}=A^{(\,^{q+1}_{\;\;n})}$$
$s_i((a_{u_1,...,u_n})_{(0\leq u_1<...<u_n\leq q-1)})
=(c_{u_1,...,u_n})_{(0\leq u_1<...<u_n\leq q)}$
where
$$c_{u_1,...,u_n}= \left\{
\begin{array}{c l}
1\; &{\rm if}\;  u_n=i\\
a_{s^i(u_1),...,s^i(u_n)}\; &{\rm if}\; u_n\neq i
\end{array}
\right.
$$
with the convention that if two consecutive indices $s^i(u_j)$, $s^i(u_{j+1})$ are equal then the corresponding element is trivial $a_{s^i(u_1),...,s^i(u_n)}=1$. We have:
\begin{theorem} $(K(n),s_i,\partial_i))$ is a simplicial group $K(A,n)$.
\end{theorem}

\begin{remark} When $n=1$ we get the classical construction of $K(A,1)$. When $n=2$ or $n=3$ we obtain the explicit description given in the previous section.
\end{remark}
\begin{remark} To understand the definition of $\partial_i$ we should remember that the elements from  $K(n)_q=A^{(\,^q_n)}$ are indexed by the n-simplices  $[u_1,u_2,...,u_n,u_n+1]$ from $\Delta_q$. Also
$\partial_i$ corresponds to $d^i:\overline{q-1}\to \overline{q}$. This means that the color of $[u_1,u_2,...,u_n,u_n+1]$ in $\Delta_{q-1}$ is the color of $[d^i(u_1),d^i(u_2),...,d^i(u_n),d^i(u_n+1)]$ from $\Delta_{q}$. If $u_n\neq i-1$ then $d^i(u_n+1)=d^i(u_n)+1$ and so the color of $[u_1,u_2,...,u_n,u_n+1]$ in $\Delta_{q-1}$ is $a_{d^i(u_1),...,d^i(u_n)}$. Otherwise we have $d^i(u_n+1)=d^i(u_n)+2$ and the color is determined by "homotopy" as described above. A similar argument can be made for the definition of $s_i$.
\end{remark}

\begin{remark} There is an obvious connection between the simplicial group $K(G,1)$ and the group cohomology $H^n(G,A)$. More exactly, $H^n(G,A)$ is the homology of the complex obtained by applying the functor $Map(\_, A)$ to the complex  associated to $K(G,1)$. The same  statement is true for the simplicial group $K(A,2)$ and the secondary cohomology $_2\,H^n(A,B)$. Similarly we can define   the ternary cohomology $\,_3H^n(B,C)$. Then for a topological space  with $\pi_1(X)=1$ and $\pi_2(X)=1$ one can construct a homotopy invariant $\,_3\kappa^5\in \,_3H^5(\pi_3(X),\pi_4(X))$.

For the general case we have to start with a 3-cocycle  $\kappa\in H^3(G,A)$, take a $\kappa$-twisted product of $K(G,1)$ and $K(A,2)$ and obtain a complex $K(G,A,\kappa^3)$. Then the secondary cohomology $\,_2H^n(G,A,\kappa;B)$ introduced in \cite{sm1} is the homology of the complex $Map(K(G,A,\kappa),B)$. At the next step, start with a 4-cocycle $\lambda \in \,_2H^n(G,A,\kappa;B)$, take a $\lambda$-twisted product between $K(G,A,\kappa^3)$ and $K(B,3)$ to obtain a complex   $K(G,A,\kappa,B,\lambda)$. Then the ternary cohomology  will be the homology of the complex $Map(K(G,A,\kappa,B,\lambda),C)$. One is then able to associate to a space $X$ an invariant $\,_3\kappa^5$ in the ternary cohomology group, and so on. This is very similar with the idea used in \cite{may} to classify simplicial groups. The main novelty is that we have an explicit way to associate to a topological space an invariant in the appropriate cohomology theory.
\end{remark}

\section{Secondary homology for commutative Hopf Algebras}

In this section $H$ is a commutative Hopf algebra. We want to associate to $H$ a cyclic object $\;_2K(H)$.  If $H$ is the group algebra $k[A]$ associated to a commutative group $A$, then $\;_2K(H)$ is the linearization of the simplicial group $K(A,2)$ described above.

Define $\;_2K(H)_q=H^{\otimes {\frac{q(q-1)}{2}}}$. An element of
$\;_2K(H)_q$ is a tensor:
$$(\otimes h_{u,v})=h_{0,1}\otimes (h_{0,2}\otimes h_{1,2})\otimes (h_{0,3}\otimes ...\otimes h_{2,3})\otimes ...
\otimes (h_{0,q-1}\otimes ...\otimes h_{q-2,q-1})$$
We define the maps $\partial_i:K_q\to K_{q-1}$ for all $0\leq i\leq q$ as follows:
\begin{eqnarray*}
\lefteqn{\partial_0((\otimes h_{u,v})_{0\leq u<v\leq q-1})=
    \varepsilon(h_{0,1}...h_{0,q-1}) h_{1,2}
 \otimes (h_{1,3}\otimes h_{2,3})}\\
&&\otimes (h_{1,4}\otimes h_{2,4}
 \otimes h_{3,4})
\otimes \; ...
\otimes (h_{1,q-1}\otimes h_{2,q-1}\otimes ...\otimes h_{q-2,q-1})
\end{eqnarray*}
\begin{eqnarray*}
\lefteqn{\partial_1((\otimes h_{u,v})_{0\leq u<v\leq q-1})=
    \varepsilon(h_{1,2}...h_{1,q-1}) \varepsilon(h_{0,1})h_{0,2}
  \otimes (h_{0,3}\otimes h_{2,3})} \\
&& \otimes (h_{0,4}\otimes h_{2,4}  \otimes h_{3,4})\otimes \; ...
\otimes (h_{0,q-1}\otimes h_{2,q-1}\otimes ...\otimes h_{q-2,q-1})
\end{eqnarray*}
for $2\leq k\leq q-1$ we define
\begin{eqnarray*}
\lefteqn{\partial_k((\otimes h_{u,v})_{0\leq u<v\leq q-1})=
    \varepsilon(h_{k,k+1}h_{k,k+2}...h_{k,q-1}) h_{0,1}
  \otimes (h_{0,2}\otimes h_{1,2}) \; ...}\\
&& \otimes (h_{0,k-1}h_{0,k}S(h_{k-1,k}^{<1>})\otimes ...\otimes h_{k-2,k-1}h_{k-2,k} S(h_{k-1,k}^{<k-1>}))   \otimes  ...\\
 &&\otimes (h_{0,q-1}\otimes h_{1,q-1}\otimes ...\otimes h_{k-1,q-1}\otimes h_{k+1,q-1}\otimes...\otimes h_{q-2,q-1})
\end{eqnarray*}
and
\begin{eqnarray*}
\lefteqn{\partial_q((\otimes h_{u,v})_{0\leq u<v\leq q-1})=
    \varepsilon(h_{0,q-1}h_{1,q-1}...h_{q-2,q-1}) h_{0,1}   \otimes (h_{0,2}\otimes h_{1,2})}\\
&&    \otimes (h_{0,3}\otimes h_{1,3}\otimes h_{2,3}) \otimes \; ...\otimes (h_{0,q-2}\otimes h_{1,q-2}\otimes ...\otimes h_{q-3,q-2})
\end{eqnarray*}

Next we define $s_i:K_q\to K_{q+1}$ for all $0\leq i\leq q$.
\begin{eqnarray*}
\lefteqn{s_0((\otimes h_{u,v})_{0\leq u<v\leq q-1})=
    1\otimes (h_{0,1}^{<1>}\otimes h_{0,1}^{<2>})}\\&&
    \otimes (h_{0,2}^{<1>}\otimes h_{0,2}^{<2>}\otimes h_{1,2}) \otimes \; ...(h_{0,q-1}^{<1>}\otimes h_{0,q-1}^{<2>}\otimes h_{1,q-1}\otimes ...\otimes h_{q-2,q-1})
\end{eqnarray*}
\begin{eqnarray*}
\lefteqn{s_1((\otimes h_{u,v})_{0\leq u<v\leq q-1})=
    1\otimes (h_{0,1}\otimes 1)}\\&&
    \otimes (h_{0,2}\otimes h_{1,2}^{<1>}\otimes h_{1,2}^{<2>})) \otimes \; ...(h_{0,q-1}\otimes h_{1,q-1}^{<1>}\otimes h_{1,q-1}^{<2>}\otimes ...\otimes h_{q-2,q-1})
\end{eqnarray*}
for $2\leq k\leq q-1$
\begin{eqnarray*}
\lefteqn{s_k((\otimes h_{u,v})_{0\leq u<v\leq q-1})=
    h_{0,1}\otimes (h_{0,2}\otimes h_{1,2})\otimes...}\\
    &&\otimes (h_{0,k-1}\otimes h_{1,k-1}...\otimes h_{k-2,k-1}) \otimes (1\otimes...\otimes 1)\\
    &&\otimes (h_{0,k}\otimes ...\otimes h_{k-1,k}\otimes 1)
    \otimes (h_{0,k+1}\otimes ...\otimes h_{k-1,k+1}
   \otimes h_{k,k+1}^{<1>}\otimes h_{k,k+1}^{<2>})\\
   &&\otimes ...
\otimes (h_{0,q-1}\otimes ...h_{k-1,q-1}\otimes  h_{k,q-1}^{<1>}\otimes h_{k,q-1}^{<2>} ...\otimes h_{q-2,q-1})
\end{eqnarray*}
and
\begin{eqnarray*}
\lefteqn{s_q((\otimes h_{u,v})_{0\leq u<v\leq q-1})=
    h_{0,1}\otimes (h_{0,2}\otimes h_{1,2})}\\&&
     \otimes \; ...\otimes (h_{0,q-1}\otimes h_{1,q-1}\otimes ...\otimes h_{q-2,q-1})\otimes (1\otimes ...\otimes 1)
\end{eqnarray*}
Finally the cyclic action is given by $\tau_q:H^{\otimes {\frac{q(q-1)}{2}}}\to H^{\otimes {\frac{q(q-1)}{2}}}$
\begin{eqnarray*}
\lefteqn{\tau_q((\otimes h_{u,v})_{0\leq u<v\leq q-1})=}\\
&&\otimes (h_{0,1}^{<1>}h_{0,2}^{<1>}...h_{0,q-2}^{<1>}h_{0,q-1}S(h_{1,2}^{<1>}h_{1,3}^{<1>}...h_{1,q-1}^{<1>}))\\
&&\otimes(h_{1,2}^{<2>}h_{1,3}^{<2>}...h_{1,q-1}^{<2>}S(h_{2,3}^{<1>}...h_{2,q-1}^{<1>})\otimes h_{0,1}^{<2>})\\
&&\otimes(h_{2,3}^{<2>}...h_{2,q-1}^{<2>}S(h_{3,4}^{<1>}...h_{3,q-1}^{<1>})\otimes h_{0,2}^{<2>}\otimes h_{1,2}^{<3>})\\
&&\otimes(h_{3,4}^{<2>}...h_{3,q-1}^{<2>}S(h_{4,5}^{<1>}...h_{4,q-1}^{<1>})\otimes h_{0,3}^{<2>}\otimes h_{1,3}^{<3>}\otimes h_{2,3}^{<3>})\\
&&...\\
&&\otimes(h_{q-3,q-2}^{<2>}h_{q-3,q-1}^{<2>}S(h_{q-2,q-1}^{<1>})\otimes h_{0,q-3}^{<2>}\otimes h_{1,q-3}^{<3>}\otimes ...
\otimes  h_{q-4,q-3}^{<3>})\\
&&\otimes (h_{q-2,q-1}^{<2>}\otimes h_{0,q-2}^{<2>}\otimes h_{1,q-2}^{<3>}\otimes ...
\otimes  h_{q-3,q-2}^{<3>})
\end{eqnarray*}

\begin{theorem} $(\;_2K(H)_q,\partial_k,s_k,\tau_q)$ is a cyclic module.
\end{theorem}
\begin{remark} If one thinks about the cyclic module $H^{(\varepsilon, 1)}$ as the cyclic module that corresponds to the first level of a "Postnikov tower", then  the second part  of that "Postnikov tower" would be a twisted product between $H^{(\varepsilon, 1)}$ and  $\,_2K(L)$. The analogy here is that $H$ plays the role of $\pi_1$ and $L$ plays the role of $\pi_2$ (therefore the need for $L$ to be commutative).
\end{remark}
\begin{remark} In order to generalize the results from \cite{cm} and \cite{c} one first needs to associate to a commutative algebra $A$ a secondary-cyclic-cohomology.  The main problem is to define an analog of the bar resolution. We will  approach that problem in a sequel to this paper.
\end{remark}



\bibliographystyle{amsalpha}

\end{document}